\mag=1200
\input amstex
\documentstyle{amsppt}
\NoBlackBoxes
\TagsOnRight
\tolerance=2000
\vsize=19cm
\topmatter
\date 11.II.2000\enddate
\address
The Lomonosov Moscow
\newline\indent
State University
\endaddress
\author
V.\,M.~Buchstaber, K.\,E.~Feldman
\endauthor
\title
An Index of an Equivariant Vector Field and Addition Theorems for
Pontrjagin Characteristic Classes
\endtitle
\abstract

The theory of indices of Morse--Bott vector fields on a manifold is
constructed and the famous localization problem for the transfer map is
solved on its base in the present paper. As a consequence, we obtained
addition theorems for the universal Pontrjagin characteristic classes
in cobordisms. These results gave us a possibility to complete
the construction, which was begun more than twenty years ago,
of the universal characteristic classes' theory.

References: 20 items.
\endabstract
\endtopmatter
\rightheadtext{An Index of an Equivariant Vector Field\dots}

\define\sk{\text{{\rm sk}}}
\define\GL{\text{{\rm GL}}}

\define\Ima{\operatorname{Im}}

\define\Tors{\operatorname{Tors}}
\define\Ind{\operatorname{Ind}}
\define\Int{\operatorname{Int}}
\define\Ker{\operatorname{Ker}}
\define\Map{\operatorname{Map}}
\define\codim{\operatorname{codim}}
\define\dima{\operatorname{dim}}
\let\Bad\relax

\let\leq\leqslant
\let\geq\geqslant

\document

\footnote""{The research was supported by RFBR, Ref. No. \,99-01-00090.
To appear in Izv. Rus. Akad. Nauk, 2000, V. 64, N 3.}

\head
\S\,1. Introduction
\endhead

The transfer map construction for fibre bundles~\cite{2},~\cite{3}
is one of the most fundamental notion in current algebraic topology.
This construction makes it possible to vitally develop and extend
applications
of the direct image method connected with such notions as the direct
image for representations, a trace of an algebraic extension, the
transfer map for vector bundles over covering, the Gysin map etc.
(see~\cite{1}--\cite{6}, \cite{10}--\cite{12},~\cite{17}).

The theory of indices of Morse--Bott vector fields on a manifold is
constructed and the famous localization problem for the transfer map is
solved on its base in the present paper. This result is remote
generalization of the classical Poincare--Hopf theorem expressing
the Euler characteristic of a manifold as the sum of zero indices
of a regular tangent vector field.

As a consequence, we obtained the formulae expressing the universal
(that is, with its values in cobordism theories) Pontrjagin characteristic
classes
of the sum of two vector bundles in terms of summands' characteristic classes
(the addition theorems). These results gave us a possibility to complete
the construction, which was begun more than twenty years ago~\cite{6},
of the universal characteristic classes' theory.

The contents of the article is as follows. In~\S\,2 we recall the
transfer map construction for fibre bundles and the related
definitions of the fiber bundle index. In~\S\,3 we introduce the
index $Ind_s(F_l)$ of a Morse--Bott vector field~$s$
on a smooth manifold~$F$ in a neighborhood of its zero submanifold~$F_l$.
This index takes its values in the cohomotopy group $\pi^0_S(F^+_l)$. We give an
example of the Morse--Bott vector field on the projective plane.
In a neighborhood of the projective line the index of this field
equals to $-1+u$, where~$u$~ is the
generator of the group $\pi^S_1=Z_2$. We prove important Theorem~3.1
of the transfer map
localization for a smooth manifold (as a fibre bundle over a point)
using the Morse--Bott vector field on the manifold. As a consequence,
we obtain an expression for the Euler characteristic $\chi(F)$ in
terms of the Euler characteristic $\chi(F_l)$ and indices $Ind_s(F_l)$
for the Morse--Bott vector field~$s$ with the collection of the connected
zero submanifolds $\{F_l\}$ (Theorem~3.2). In~\S\,4 we prove
in general case the localization theorem for the transfer map of an arbitrary
fiber bundle with a smooth fiber. We demonstrate that the known
results~\cite{5},~\cite{11} in this direction can be reduced to special
cases of the theorem.

The rest part of the work is devoted to the applications of the results
obtained to the theory of Pontrjagin characteristic classes of real vector
bundles. In~\S\,5 we construct an important Morse--Bott vector field over
grassmannization of the Whitney sum of two real vector bundles.
We describe connected components of the zero submanifolds and calculate
the indices of this field in their neighborhoods.

In~\S\,6 we obtain one of the central results of the article:
the addition theorem for the universal Pontrjagin classes of real vector
bundles. In details we investigate the case of the first half--integer
class~$p_{1/2}$. We claim that on the image of the map
$$
p_{1/2}\:KO^*(X)\to U^2(X)
$$
there is a new algebraic operation that cannot be reduced to any formal
group. Call attention to the fact that in our short paper~\cite{9} we made
an inaccuracy in the formula for $p_{1/2}(\xi\oplus\zeta)$ in complex
cobordisms, which we correct in the present work.

In the last section we investigate Pontrjagin characteristic classes of
stable complex vector bundles. We prove that if a real vector bundle
admits a stable complex structure, its half--integer Pontrjagin classes are
equal to zero while its integer Pontrjagin classes have a realization
in symplectic cobordisms.

The authors express their deep gratitude to Dr. I.\,A.~Dynnikov for
helpful discussions and valuable suggestions during this work was being
executed.

\head
\S\,2. The Index of a Smooth Fibering
\endhead

In this section we shall give definitions of the main notions connected
with the transfer map of Becker--Gottlieb~\cite{3}.

Let us consider two topological spaces~$X$ and~$Y$ with marked points and
the set of homotopic classes $[X,Y]$ of all maps from~$X$ to~$Y$,
sending a marked point to another marked one. We denote by $\{X,Y\}$ the group
of stable homotopic classes of maps
$$
\{X,Y\}=\lim_{n\to\infty}[S^n\wedge X,S^n\wedge Y].
$$

\definition{Definition 2.1}
{\it A zero--dimensional cohomotopy group of a space~$X$
with a marked point\/} is the group
$$
\pi^0_S(X)=\{X,S^0\}=\lim_{n\to\infty}[S^n\wedge X,S^n].
$$
\enddefinition

Let $X^+$ denote a disjoint union of the space~$X$ and a point,
which we shall assume to be a marked point. If the marked point in the
space~$X$ is already given, the space $S^n\wedge X^+$ is homotopy--equivalent
to $S^n X\vee S^n$, and therefore,
$$
\pi^0_S(X^+)=\pi^0_S(X)\oplus Z.
$$
We shall denote by~$\epsilon$ the canonical projection $\pi^0_S(X^+)\to Z$.
If the space~$X$ is the sphere $S^k$,
$$
\Ker(\epsilon) \cong \pi^S_k
$$
where $\pi^S_k$ is the $k$-th stable homotopy group of spheres.

Every fibering $(E,F,B,p)$ with the smooth compact closed fiber~$F$
defines the canonical element in the group $\{B^+,E^+\}$. Existence of this
element
is provided by the transfer map of Becker and Gottleib~\cite{3}.
Let us give the construction of this element in case of a smooth fibering~$E$ with the compact base (for arbitrary cell bases see the definition in
~\cite{8}).

The space of the fibering~$E$ can be embedded into $\Bbb R^n$ for
suitable~$n$. Let $i$ denote this embedding. We consider the space $\tau_F(E)$ of tangents
along the fiber of the bundle~$E$  and the normal bundle $\nu(E)$ of the
inclusion $p\times i\: E\subset B\times\Bbb R^n$. It is easy to check that
$\tau_F(E)\oplus \nu(E)$ is a trivial $n$-dimensional vector bundle
over~$E$~\cite{3}. Let $T(\xi)$ denote the Thom space of an arbitrary
vector bundle~$\xi$.

\definition{Definition 2.2}
{\it The transfer map of a smooth fiber bundle $(E,F,B,p)$\/}
is the composition of the maps
$$
\tau(p)\:S^n\wedge B^+\to T\nu(E)\to
T\bigl(\nu(E)\oplus\tau_F(E)\bigr)=T[\bar n]
$$
where the first map is the Pontrjagin--Thom map representing the glueing
of the tubular neighborhood's exterior of the submanifold~$E$ in the manifold
$B\times\Bbb R^n$, while the second map is the inclusion of the fiber bundle
$\nu(E)$ on the direct summand.
\enddefinition

Notice that the Thom space of the trivial~$n$-dimensional vector
bundle over~$B$ is the~$n$-fold suspension over~$B^+$. Therefore, the transfer map
defines a certain element of the group $[S^n\wedge B^+,S^n\wedge E^+]$.
The image of this element in the group $\{B^+,E^+\}$ is independent of a
choice of the embedding~\cite{3}. Thus, the transfer map
correctly defines an element, which we shall
denote by $\bigl\{\tau(p)\bigr\}$, of the group $\{B^+,E^+\}$.

For any cohomology theory $h^*(\cdot)$ a stable map from~$X$ to~$Y$
uniquely defines groups' homomorphism
$$
h^*(Y)\to h^*(X).
$$
We shall denote by $\tau(p)^*$ the corresponding homomorphism for the transfer
map of the fibering $(E,F,B,p)$. Let $\pi$ represent the projection on
the first summand $\pi:S^n\wedge E^+\to S^n$.

\definition{Definition 2.3}
{\it An index $I(E)$ of the fiber bundle $p:E\to B$\/} is an
element $\bigl\{\pi\circ\tau(p)\bigr\}\in\pi^0_S(B^+)$.
An index of this fiber bundle
in the cohomology theory $h^*(\cdot)$ is the element
$I_h(E)=\tau(p)^*(1)\in h^0(B^+)$ (we consider the subgroup~$h^0$ in the
group~$h^*$ as a ring with unit). Thus, $I(E)$ is the fiber bundle's index
in the stable cohomotopy theory  $\pi^*_S(\cdot)$.
\enddefinition

 Notice that for any cohomology theory $h^*(\cdot)$ with unit
there is a canonical transformation of cohomology theories
(see~\cite{19})
$$
\mu_h: \pi^*_S(\cdot)\to h^*(\cdot)
$$
which transforms the unit of the theory $\pi^*_S(\cdot)$ into the unit of
the theory $h^*(\cdot)$. In the theory $h^*(\cdot)$ the index of the fiber
bundle $ p:E\to B$  is equal to $\mu_h\bigl(I(E)\bigr)$.

The indices of the fiber bundles are already non--trivial in the simplest
case for the fiberings over a point. In that case $\pi^0_S(B^+)\cong Z$,
and the index of the fibering is equal to the Euler characteristic of the
fiber~\cite{3}.

\head
\S\,3. Localization
\endhead

The index of the tangent vector field on the smooth manifold in a
neighborhood of a nondegenerate zero is the fundamental characteristic
of the zero considered~\cite{14}. This index is defined as
determinant's sign of the linear part of the vector field in a small
neighborhood of the zero~\cite{14}. In this section we claim
that in a neighborhood of the degenerate zero the differential of the field
gives an important characteristic of this zero. In terms of the
characteristics introduced one can express the coefficients, for which the
transfer map
of the smooth manifold can be represented by the transfer maps' linear
combination of the zero submanifolds of the vector field on the initial
manifold.

 Before proving the main theorem we shall give an example that will help one
to understand the transfer map localization property.
Let us begin with necessary notation.

Each element of the group $\pi^0_S(B^+)$ can be determined by a continuous map
as $S^n\wedge B^+\to S^n$ for~$n$ large enough. Let us denote by~$\pi$
the projection on the first factor in the space $S^n\wedge B^+$.
The representative~$\gamma$ of an element of the group $\pi^0_S(B^+)$
defines a map
$$
f_{\gamma}\:B\to\Map(S^n,S^n)
$$
by the formula $f_{\gamma}(b)(y)=\gamma(y,b)$. The inverse is also true. Any
continuous map $f:B\to Map(S^n,S^n)$ defines a certain element of
$\pi^0_S(B^+)$ by the formula
$\gamma_f=\bigl\{\pi\bigl(f(b)(y),b\bigr)\bigr\}$. Let us
identify~$S^n$ with $\Bbb R^n\cup\{\infty\}$ and define the action
$\GL(n, \Bbb R)$ on~$S^n$ by this identification. Even in case of
$f(b)\in\GL(n,\Bbb R)$ for any
$b\in B$ the correspondent element of the cohomotopy group can take
nonzero value after the projection onto $\Ker(\epsilon)$.

\example{Example 3.1}
Let $B=S^1$ and $\xi_1$ denote the nontrivial
one--dimensional real vector bundle over~$S^1$. We have
$\xi_1\oplus \xi_1$ is the trivial two-dimensional real vector bundle.
Let us consider the map
$$
S^2\wedge S^{1+}=T(\xi_1\oplus\xi_1)\overset{(1,-1)}\to{\longrightarrow}
T(\xi_1\oplus\xi_1)=S^2\wedge S^{1+}\overset{\pi}\to{\longrightarrow}S^2.
\tag 3.1
$$
Since $S^2\wedge S^{1+}\sim S^3\vee S^2$, the map constructed above
defines two maps $\alpha_1:S^3\to S^2$ and $\alpha_2:S^2\to S^2$.
It is easy to check that~$\alpha_2$ is the reflection relatively to the
central
plane. Let us show that~$\alpha_1$ defines a generator of the group
$\pi_3(S^2)$. Identify~$S^2$ and $\Bbb C^1\cup \{\infty\}$ setting the infinite
point to be a marked one. We consider the preimage of the zero under the
map~$\alpha_1$.
It is the zero section in the bundle $\xi_1\oplus\xi_1$, that is, the circle.
The whole preimage $\Bbb C^1\subset S^2$ under the map $\alpha_1$ is
the space of the vector bundle $\xi_1\oplus\xi_1$. The splitting
of the trivial two-dimensional real vector bundle over~$S^1$ in the direct
sum $\xi_1\oplus\xi_1$ defines the chosen basis in each fiber of the
bundle
$$
\bigl\langle e_1(\phi),e_2(\phi)\bigr\rangle=
\langle e^{i\phi/2},ie^{i\phi/2}\rangle
$$
where $\phi\in[0,2\pi]$ is the angular coordinate on the zero section
(that is, on~$S^1$). For any point $\{\phi,x+iy\}$ from
$S^1\times \Bbb C^1\subset S^3$ its image under the map~$\alpha_1$ is
$$
\alpha_1\bigl(\{\phi,x+iy\}\bigr)=\bigl\{(x-iy)e^{i\phi}\bigr\}\in C^1.
$$
In coordinates $(z_1,z_2)\in \Bbb C^1\times \Bbb C^1$,
\ $|z_1|=1$, the map $\alpha_1$ can be expressed by the formula
$$
\alpha_1(z_1,z_2)=\bar z_2z_1,
$$
which shows that the map $\alpha_1$ is homotopic to the projection in
the Hopf bundle.

Therefore, the element of the group $\pi^0_S(S^{1+})$ corresponding to the
map~\thetag{3.1} is $-1+u$, where~$u$ is the generator of the group
$\pi^S_1\cong Z_2$.
\endexample

Now we formulate the localization property, which at first
will be proved for smooth manifolds, that is, in case of fiber bundles over a
point.

We shall use the notion of a tangent vector field of the Morse--Bott type
(see~\cite{16}); Morse vector fields are the special case of such
fields.

\definition{Definition 3.1}
A tangent vector field on a smooth manifold~$F$
is called {\it a Morse--Bott field} in case of the two following conditions
are satisfied:

1) the zero set of the vector field forms a finite set of connected compact
closed submanifolds in~$F$;

2) the kernel of the Jacobian matrix of the vector field at the points of the
zero submanifold coincides with the tangent space of the zero submanifold.
\enddefinition

Let us consider the connected component~$F_l$ of the zero set of the
Morse--Bott
vector field. Let~$\nu_l$ be the normal bundle of the embedding
$F_l\subset F$. The bundle space~$\nu_l$ can be identified
with the tubular neighborhood of the submanifold~$F_l$ in~$F$ using
the exponential map~\cite{14}. The restriction of the Morse--Bott vector
field on the tubular neighborhood defines the vector field on the
space~$\nu_l$. Owing to the condition~2) the last vector field must have everywhere
non-zero projection on the bundle of tangents along the fibers of
the bundle~$\nu_l$ in a small enough neighborhood of the zero section.
Thus, without loss of generality, it can be assumed that the vector field
constructed on the space~$\nu_l$ is tangent along the fibers of the
bundle~$\nu_l$ everywhere. Let us denote by~$s_l$ this vector field
$$
s_l\:\nu_l\to\hat\tau(\nu_l)
$$
where $\hat \tau(\nu_l)$ is a bundle of tangents along the fibers. Let~$p$ be
the projection in the bundle~$\nu_l$, then~$\hat\tau(\nu_l)\cong p^*\nu_l$
(see, for example,~\cite{13}). Thus, the vector field~$s_l$ allows
us to construct the map $\psi^s_l\:\nu_l\to \nu_l$ by the rule
$$
\psi^s_l(v)=\hat ps(v)
$$
where $\hat p$ is the bundle map $\hat p\: p^*\nu_l\to \nu_l$.

 Let us consider a smooth compact closed manifold~$F$ and an arbitrary
Morse--Bott vector field~$s$ on~$F$. Let $F_1,\dots,F_m$ be the
whole collection of connected compact closed zero submanifolds of the
vector field~$s$ on~$F$. Let us fix the embedding of the manifold~$F$
in the~$n$-dimensional Euclidean space:
$$
F_1,\dots,F_m\subset F\subset\Bbb R^n.
$$

Using the vector field~$s$ we correspond  the
map
$$
j_l\: S^n\wedge F^+_l\to S^n\wedge F^+_l
$$
to each submanifold~$F_l$ as follows. The embedding of the smooth manifolds
$ F_l\subset F\subset \Bbb R^n $ induces the splitting of the
trivial~$n$-dimensional vector bundle over~$F_l$ into the direct sum of
the tangent bundle~$\tau_l$ of~$F_l$, the normal bundle~$\nu_l$ of
the embedding $F_l\subset F$, and the restriction on~$F_l$ of the
normal bundle~$\nu$ of the embedding~$F\subset\Bbb R^n$. The field~$s$
determines the map $\phi^s_l\: \nu_l\to\nu_l$. Let us define the map~$j_l$
on the subbundles~$\nu$ and~$\tau_l$ by identical, Ánd that on~$\nu_l$
by~$\phi^s_l$.

By the above, let $\pi$ be the projection $S^n\wedge F^+_l\to S^n$
on the first factor.

\definition{Definition 3.2}
An element $\{\pi\circ j_l\}\in\pi^0_S(F^+_l)$ will be called {\it an index
$\Ind_s(F_l)$ of the vector field~$s$\/} in a
neighborhood of the zero submanifold~$F_l$.
\enddefinition

\remark{Remark 3.1}
An index of the isolated zero of the vector field
coincides with $\epsilon\bigl(\Ind_s(pt)\bigr)$
and is equal to the ordinary index
of the isolated zero~\cite{14}.
\endremark

\example{Example 3.2}
Using Example~3.1 we claim that an index
$Ind_s(F_l)\in \pi^0_S(F^+_l)$ can have non--zero component in
$\Ker(\epsilon)$.

Let us define an action of the positive real numbers group~$\Bbb R_+$ on
$RP^2$ by the formula
$$
t\circ(x_1:x_2:x_3)= (x_1:x_2:t^{-1} x_3).
$$
The differential of this action  determines the tangent vector field~$s$
on $RP^2$ if $t=1$. All fixed points of the action
are zeros of~$s$, that is,
$$
RP^1=\bigl\{(x_1:x_2:0)\bigr\}\subset RP^2,
$$
and the point $(0:0:1)\in RP^2$. The normal bundle of the embedding
$RP^1\subset RP^2$ can be identified with the nontrivial one-dimensional
real vector bundle~$\xi_1$ over~$S^1$. In a neighborhood
$RP^1=S^1\subset RP^2$ the vector field~$s$ is determined by the
formula $ s(v)=-v$, \ $v\in \xi_1\approx N(S^1)\subset RP^2$. Owing
to Example~3.1 an index $\Ind_s(RP^1)$ is equal to $-1+u$, where~$u$
is the generator of the group $\pi^S_1\cong Z_2$.
\endexample

 Let the map $\tau_l$ be the transfer for the manifold~$F_l$. Let us denote
by $i$ the Thom space embedding of the trivial~$n$-dimensional
vector bundle over $F_1\cup\dots\cup F_m$ in the Thom space of the
trivial~$n$-dimensional vector bundle over~$F$ that corresponds to the
embedding
$$
F_1\cup\dots\cup F_m\subset F.
$$

\proclaim{Theorem 3.1 \rm( the localization property)}
The transfer map~$\tau$ for the manifold~$F$ is homotopic to the composition
$$
i\circ(j_1\circ\tau_1\vee\dots\vee j_m\circ\tau_m)\circ g
$$
where~$g$ is the map sending the sphere~$S^n$ into the union of~$m$ spheres.
\endproclaim

\demo{Proof}
The embedding of the Thom space
$$
\alpha\:T\nu(F)\to T\bigl(\nu(F)\oplus\tau(F)\bigr)
$$
from the definition of the transfer map can be changed for a
homotopic one using the tangent field. In fact, let us consider
a one--parametric family of maps
$$
\alpha_t(v)=\alpha(v)\oplus ts\bigl(\pi(v)\bigr),\qquad v\in\nu(F),
$$
where $\pi$
is the projection in the bundle $\nu(F)$. Choosing a large enough
parameter~$t$ we shall get the  map, homotopic to~$\alpha$, sending
small neighborhoods' exteriors of the zeros of the field~$s$
at the marked point of the space $T\bigl(\nu(F)\oplus\tau(F)\bigr)$.
This means the map~$\tau$ can be passed through the map~$g$ and the union of
maps sending~$S^n$ to $T\nu(F_l)$, \ $l=1,\dots,m$.

To describe properties of the map $\tau$ in a neighborhood of the
zeros of the field $s$ it will be conveniet to change the map $\alpha$
on a homotopic one using the canonical identification of tubular
neighborhoods of the
submanifolds $F_l$  and normal bundles of embeddings
$F_l\subset F$, \ $l=1,\dots,m$.

We denote by $exp_l$ the exponential map that defines a diffeomorphism
of unit-disk bundle associated with the normal bundle $\nu_l$ of
the embedding $F_l\subset F$ and the tubular neighborhood
$N(F_l)\subset F$ of the submanifold~$F_l$ in~$F$, $l=1,\dots,m$.
 We define two functions
$$
\lambda_l,\mu_l\:D^{k_l}\to \Bbb R^{k_l}
$$
for each $l=1,\dots,m$,
where~$D^s$ denotes the unit disk of dimension~$s$, \ $k_l=\codim F_l$,
setting
$$
\lambda_l(v)=
\cases
|v|,&\quad|v|\leq\frac12,
\\
1-|v|,&\quad|v|>\frac12
\endcases
$$
and
$$
\mu_l(v)=
\cases
0,&\quad|v|\leq\frac12,
\\
2|v|-1,&\quad|v|>\frac12.
\endcases
$$
Let us consider the map
$$
\gather
\tilde\alpha\:F\to\tau(F),
\\
\tilde\alpha(x)=
\cases
\alpha(x),&\quad x\notin N(F_l),
\\
\bigl\{\exp\bigl(\mu_l\bigl(\exp^{-1}_l(x)\bigr)\exp^{-1}_l(x)\bigr),
\lambda_l\bigl(\exp^{-1}_l(x)\bigr)\exp^{-1}_l(x)\bigr\},&\quad
x\in N(F_l),
\endcases
\\
l=1,\dots,m.
\endgather
$$
Since $\lambda_l$ is homotopic to the zero map and
$\mu_l\bigl(\exp^{-1}_l(x)\bigr)\exp^{-1}_l(x)$ is homotopic to the map
$\exp^{-1}(x)$, the map~$\tilde \alpha$ is homotopic to the embedding on
the zero section $F\to \tau(F)$.

We remark that
$$
\nu(F)\bigr|_{N(F_l)}\cong\psi_l(1)^*\nu(F)\bigr|_{F_l}
$$
where $\psi_l(t)$ is the retraction of the normal~$t$-tube $N_t(F_l)$
on~$F_l$. Consequently, the map~$\alpha$ is homotopic to the map
$$
\beta(v)=
\cases
v,&\quad \pi(v)\notin N(F_l),
\\
\psi_l\bigl(\bigl|\pi\alpha\bigl(\pi(v)\bigr)\bigr|\bigr)^*(v),&\quad
\pi(v)\in N(F_l),
\endcases\qquad
l=1,\dots,m.
$$

 Thus, the map constructed is the result of a small embedding
perturbation
$$
T\nu(F)\to T\bigl(\nu(F)\oplus\tau(F)\bigr)
$$
under which the tubular neighborhoods' fibers of the submanifolds~$F_l$
transfer into the fibers of the normal bundles of the embeddings
$F_l\subset F$ corresponding them under the exponential map.

Let us construct the map
$$
\tau(t)=\bigl(\beta(v)\oplus t\tilde\alpha_*
\bigl(s\bigl(\pi(v)\bigr)\bigr)\bigr)
$$
where $v\in T\nu(F)$, \ $t\in [0;\infty)$, and we turn~$t$ to $\infty$.
We obtain the map $\tau(\infty)$ that sends at infinity the points of the
bundle $\nu(F)$ lying in the fibers over exterior of the tubular neighborhoods
of the submanifolds~$F_l$.

 The fiber of the tubular neighborhood of the submanifold~$F_l$ in~$F$
under $\tau(\infty)$ maps onto the fiber of the subbundle~$\nu_l$
(see above).  One can neglect the first summand in the sum
$$
\beta(v)\oplus t\tilde\alpha_*\bigl(s\bigl(\pi(v)\bigr)\bigr)
$$
for large~$t$.
Thus, the restriction of the map $\tau(\infty)$ to the fiber of
the tubular neighborhood of the submanifold~$F_l$ in~$F$ coincides
with the map $\phi^s_l$ (see above). To complete the proof we have
to remark that the image of the map $\tau(\infty)$ coincides with the
image of the map~$i$.
\enddemo

\remark{Remark 3.2}
The conditions imposed on the vector fields in the
formulation of the localization property can be weakened. It is enough
to require the vector field to be homotopic to the Morse--Bott vector
field in the class of the vector fields with the fixed zero set.
\endremark

\remark{Remark 3.3}
It is easy to check that if the vector field~$s$
is homotopic to the outward normal field on the boundary of the tubular
neighborhood of the zero submanifold~$F_l$, the map $j_l$ corresponding
to that submanifold is homotopic to the identical map~\cite {8}. In
particular, this situation always holds in case of the vector field defined
by a circle action on the manifold~\cite{2}.
\endremark

\remark{Remark 3.4}
Immediately from the proof of Theorem~3.1 it
follows that if the tangent vector field on the manifold has no zeros,
the transfer map for this manifold is homotopic to the map into
the point.
\endremark
The localization property and Remark~3.1 allow to immediately obtain
the following generalization of Poincare--Hopf Theorem (compare~\cite{20}).

\proclaim{Theorem 3.2}
Let~$F$ denote a smooth compact closed manifold.
Let the Morse--Bott vector filed~$s$ be given on~$F$. Assume that zeros
of the vector field~$s$ form a finite collection of connected compact
closed submanifolds  $F_1,\dots,F_m$ of the manifold~$F$. Then
$$
\chi(F)=\sum^m_{l=1}\epsilon\bigl(\Ind_s(F_l)\bigr)\cdot\chi(F_{l})
$$
where $\chi(F)$ is the Euler characteristic of the manifold~$F$.
\endproclaim

\demo{Proof}
We shall use the notation of Theorem~3.1. At first we shall
consider an arbitrary tangent vector field on~$F$ with isolated
nondegenerated zeros.
Let us apply the localization property to this field. We obtain
$$
I(F)=\{\pi\circ\tau\}=\sum^N_{k=1}\{\pi\circ j_k\}=
\sum^N_{k=1}\epsilon\bigl(\Ind_s(z_k)\bigr)
$$
where~$z_k$ are zeros of the vector field, $k=1,\dots,N$. Owing to
Remark~3.1 this means that $I(F)=\chi(F)$ is the Euler characteristic of the
manifold~$F$.

Applying the localization property to the vector field~$s$ we obtain
$$
\split
\chi(F)&=I(F)=\{\pi\circ\tau\}=\sum^m_{l=1}\{\pi\circ j_l\circ\tau_l\}=
\sum^m_{l=1}\tau(p_l)^*\{\pi\circ j_l\}
\\
&=
\sum^m_{l=1}\epsilon\bigl(\Ind_s(F_l)\bigr)\cdot I(F_l)=
\sum^m_{l=1}\epsilon\bigl(\Ind_s(F_l)\bigr)\cdot\chi(F_l).
\endsplit
$$
\enddemo

\head
\S\,4. Equivariant Vector Fields and the Transfer Map of Becker--Gottlieb
\endhead

In this section we shall generalize the localization property for
the transfer map to the case of an arbitrary fiber bundles with smooth
fiber. The known results \cite{5},~\cite{11} are the particular
cases of our investigation. In this connection we show the way to obtain
these results in terms of the vector field index introduced above.

We consider smooth principal~$G$-bundle $(\Cal E,G,B,p)$ whose base space
is a compact connected manifold, and  $G$ is a compact Lie group. Let
$F$ denote a compact closed connected manifold equipped with a smooth
$G$-action. We assume that there is a tangent Morse--Bott vector field~$s$,
that is invariant under $G$-action, on~$F$ . Zeros of the vector field~$s$
represent themselves the union of irreducible invariant submanifolds of
the manifold~$F$. Let us denote them by $F_1,\dots,F_m$. Because of
irreducibility the following manifolds are connected
$$
E_l=\Cal E\times_GF_l.
$$
Every invariant under $G$-action tangent Morse--Bott vector field
$s$ on~$F$ allows us to construct a vector field on $E={\Cal E}\times_GF$,
that is tangent along the fibers of the bundle~$E$. We denote it
also by~$s$. We have
$$
s\:E\to \tau_F(E)=\Cal E\times_G\tau(F).
$$
More over the vector field~$s$ on~$E$ is also a Morse--Bott field.
Zeros of the field~$s$ on~$E$ represent themselves subbundles
$E_l={\Cal E}\times_GF_l$, \ $l=1,\dots,m$, of the bundle~$E$.
Let us fix an embedding of the bundle~$E$ in~$\Bbb R^n$. Exactly
as it was in the previous section the vector field~$s$ defines the maps
$$
j_l\:S^n\wedge E^+_l\to S^n\wedge E^+_l.
$$
 Let $\tau(p_l)$ denote the transfer map for the bundle~$E_l$.
Let $\hat g$ be the map sending  $S^n\wedge B^+$ into the union of~$m$
spaces $S^n\wedge B^+$ which is induced by the map $g\:S^n\to \bigvee_{m} S^n$.

\proclaim{Theorem 4.1}
The transfer map $\tau(p)$ for the bundle~$E$ is
homotopic to the composition of the maps
$$
\tau(p)\sim i\circ\bigl(j_1\circ\tau(p_1)\vee\dots
\vee j_l\circ\tau(p_l)\bigr)\circ\hat g
$$
where~$i$ is an embedding of the Thom space of the trivial~$n$-dimensional
vector bundle over $E_1\cup\dots\cup E_m$ into the Thom space of the
trivial $n$-dimensional vector bundle over~$E$.
\endproclaim

\demo{Proof}
Since the vector field~$s$ is invariant respectively
to the action of the group~$G$ on~$F$, all the homotopies described in
the proof of Theorem~3.1 are equivariant accordingly to the action
of the structural group on the fiber. Thus, the equivariant analog of
Theorem~3.1 holds. As a matter of fact this analog is the formulation
of Theorem~4.1.
\enddemo

\remark{Remark 4.1}
The conditions on the base space of the bundle~$E$ can be
weakened. We may only require the base to be a cell complex.
In this case it is enough to carry on the proof for the fiber of the bundle
and then use its equivariantness respectively to the action
of the structural group of the bundle. Let us note that homotopy equivalence
for an infinite cell complex is understood as equivalence of the map's
stable classes.
\endremark

Let $h^*(\cdot)$ be a multiplicative cohomology theory with unit.
We denote by $\bigl(\Ind_s(E_l)\bigr)_h$ the image of $\Ind_s(E_l)$ under
the canonical transformation of cogomology theories $\mu_h\:\pi^*_S\to h^*$.

\proclaim{Theorem 4.2}
Under the conditions described above the following
formula holds
$$
\tau(p)^*a=\sum^m_{i=1}\tau(p_l)^*
\bigl(\bigl(\Ind_s(E_l)\bigr)_h\cdot i^*_la\bigr)
$$
where $a\in h^*(E^+)$ and $i_l\:E_l\subset E$ is an embedding of the bundles,
$l=1,\dots,m$.
\endproclaim

\demo{Proof}
Let us notice that the maps~$j_l$, as elements from
the group of stable-equivalent maps from~$E_l$ into itself, uniquely
determine the gomomorphism~$j^*_l$ of the cogomology groups of the
corresponding spaces. Directly from the definition of multiplication in
multiplicative cohomology theories it follows that the gomomorphism~$j^*_l$
is multiplication on $\bigl(\Ind_s(E_l)\bigr)_h$.
Thus, the statement of Theorem~4.2 follows from Theorem~4.1.
\enddemo

Theorem 4.2 allows us to obtain the result of~\cite{5} in the
following form.

\proclaim{Corollary 4.1}
For gomomorphisms in the singular cogomology theory
the following equality holds
$$
\tau(p)^*=\sum^m_{l=1}\epsilon\bigl(\Ind_s(E_l)\bigr)\tau(p_l)^*i^*
$$
where $i\:E_l\subset E$, \ $l=1,\dots,m$, are embeddings of zero subbundles
of the field~$s$.
\endproclaim

To obtain the result of~\cite{11} from Theorem~4.2
we shall start with discussion of necessary definitions.

Let us consider the principle~$G$-bundle $({\Cal E},G,B,p)$, where~$G$ is
the compact Lie group. Let~$G$ act on a smooth closed compact manifold~$F$.

\definition{Definition 4.1}
We shall say that {\it orbits of points $x_1,x_2\in F$ belong to the same type}
if their stabilizers $N(x_1), N(x_2)$ are conjugated.
The orbit type of the point~$x_1$ is less then that of the point~$x_2$ if
$N(x_2)$ is conjugated to a subgroup in the group $N(x_1)$.
\enddefinition

 The union of the points belonging to the orbits of the same type~$\gamma$
forms a smooth submanifold~$F_{\gamma}$ in~$F$ whose boundary
belongs to the orbit of a smaller type. Let us consider the factor space
$Y=F/G$. The closure of the set $Y_{\gamma}= F_{\gamma} / G \subset Y$ is
a smooth manifold with angles and its boundary lies in
$\bigcup_{\gamma'<\gamma} Y_{\gamma'}$~\cite{18}.

\proclaim{Theorem 4.3 \cite{18}}
There is a simplicial subdivision of
space~$Y$, for which the interior of each simplex lies in $Int(Y_{\gamma})$
for some~$\gamma$.
\endproclaim

Let us show how one can prove Theorem~5.14 from~\cite{11} using
Theorem~4.1.

\proclaim{Theorem 4.4}
Under the conditions described above in an arbitrary
cogomology theory for a gomomorphism induced by the transfer map
of the bundle ${\Cal E}\times_GF$, the following formula holds
$$
\tau(p)^*=\sum_{\sigma\subset Y}(-1)^{\dima\sigma}\tau(p_{\sigma})^*
$$
where $\tau(p_{\sigma})$ is the transfer map
$\Cal E\times_G\bigl(G/N(x_{\sigma})\bigr)$, \ $x_\sigma$
is an arbitrary preimage
of a simplex $\sigma\subset Y$ barycenter in the simplicial subdivision
constructed in Theorem~{\rm4.3}.
\endproclaim

\demo{Proof}
Let us construct the canonical tangent Morse--Bott vector field
on the manifold~$F$. To accomplish this at first we shall construct a Morse vector field
on each~$n$-dimensional skeleton  $\sk^n Y$ in the simplicial subdivision
of the space~$Y$. The construction will be held by induction on skeleton's
dimension. Assume that zeroes of the field~$s$ are all simplex vertices
of the simplicial subdivision. Define the vector field~$s$ on $\sk^1 Y$ to be
repelling from $\sk^0 Y$ and attracting to the barycenter of one-dimensional
simplexes. Assume that we have constructed the vector field on $\sk^n Y$.
Extend it to the field on $\sk^{n+1} Y$ supposing barycenters of
$(n+1)$-dimensional simplexes to be attracting zeroes of the field~$s$. This
completes both the induction step and the construction of the
field. Thus, we have constructed the Morse vector field on each skeleton
$\sk^n Y$, whose zeroes are all the barycenters of the simplexes. A zero
index in the barycenter of a~$n$-dimensional simplex is equal to $(-1)^n$.
The map $g\: F\to F/G=Y$ induces the trivial bundle
$g\:g^{-1}\bigl(\Int(\sigma)\bigr)\to \Int(\sigma)$ for any simplex
$\sigma\subset Y$.
Fix a Riemannian metric, invariant relatively to~$G$, on the manifold~$F$.
The canonical splitting of the tangent bundle
$$
\tau\bigl(g^{-1}\bigl(\Int(\sigma)\bigr)\bigr)=
g^*\tau\bigl(\Int(\sigma)\bigr)\oplus\Ker g_*
$$
allows one to lift the vector field~$s$ on the manifold
$g^{-1}\bigl(\Int(\sigma)\bigr)$.
This lifting will be compatible for different simplexes since
$$
\partial g^{-1}(\sigma)=g^{-1}(\partial\sigma).
$$
Let us denote by $\hat s$ the lifting of the field~$s$ on the whole~$F$.
The vector field~$\hat s$ is the Morse--Bott vector field on~$F$ that is
invariant respectively to the action of the group~$G$ on~$F$. Submanifolds
$G/N(x_{\sigma})$ are zeroes of the vector field~$\hat s$. We have
$$
\Cal E\times_Gg^{-1}(\sigma)\cong
\bigl({\Cal E}\times_GG/N(x_{\sigma})\bigr)\times\sigma.
$$
This means the normal bundle of an embedding
$$
\bigl({\Cal E}\times_GG/N(x_{\sigma})\bigr)\subset{\Cal E}\times_GF
$$
is a trivial vector bundle. Thus,
$\Ind_{\hat s}\bigl(G/N(x_{\sigma})\bigr)=(-1)^n$, where~$n$~
is the dimension of the simplex~$\sigma$. Hence, the statement of the
theorem follows.
\enddemo

\example{Example 4.1 \rm(compare~\cite{10})}
Let us consider a fiber bundle
of the Klein bottle $KL$ over the circle~$S^1$ with fiber~$S^1$
$$
p\:KL\overset{S^1}\to{\longrightarrow}S^1.
$$
The index $I(p)$ of the bundle constructed is equal to the generator of the
group
$$
Z_2\cong\pi^S_1\subset\pi^0_S(S^{1+}).
$$
\endexample

\demo{Proof}
Let $\xi_1$ be a nontrivial one-dimensional real vector bundle
over~$S^1$. Then the bundle $p\:KL\to S^1$ is isomorphic to the bundle
$p\: RP(\xi_1\oplus 1)\to S^1$. The structural group of this bundle acts
on the fiber~$S^1$ by reflections respectively to the line fixed.
The fixed points of the action correspond to two subbundles $RP(\xi_1)$ and
$ RP(1)$ in the bundle $ RP(\xi_1\oplus 1)$. The normal bundles of embeddings
$ RP(\xi_1)$, \ $ RP(1)\subset  RP(\xi_1\oplus 1)$ are isomorphic to the
bundle $\xi_1\otimes 1\cong \xi_1$. Let us construct the vector field on~$S^1$
that is invariant respectively to the action of the structural group of
the bundle $ RP(\xi_1\oplus 1)$. The fixed points of the action are zeroes of
this vector field. The field is repelling from one zero
and attracting to another one. Besides, the vector field is
symmetric relatively to the line passing through the fixed points of the action.
This field defines the fiber Morse--Bott vector field~$s$ on
$ RP(\xi_1\oplus 1)$. Zeroes of this vector field are the subbundles
$ RP(\xi_1)$ and $ RP(1)$, whose bundle spaces are the sections of
the bundle $ RP(\xi_1\oplus 1)$, that is,~$S^1$. The localization property
(Theorem~4.1) applied to the vector field~$s$ gives
$$
I(p)=\Ind_s\bigl( RP(\xi_1)\bigr)+\Ind_s\bigl( RP(1)\bigr).
$$
The repelling zero of the field~$s$ has the index that is equal to~1. The
attracting zero has the index as in Example~3.2, that is, $-1+u$, where~$u$ is
the generator of the group~$\pi^S_1$. Thus, $I(p)=u$.
\enddemo

\head
\S\,5. Morse--Bott Vector Fields on Grassmannizations of Splittable
Vector Bundles
\endhead

 In this section we shall construct the fiber Morse--Bott vector field on
grassmannization of the splittable vector bundle and shall explicitly
calculate its indices. Properties of this vector field play a significant
role while proving addition formulae for Pontrjagin characteristic classes.

 Below we shall denote by $RG_k^n$ the Grassmannian manifold of
$k$-dimensional subspaces of a $n$-dimensional real vector space.
Consider a $n$-dimensional real vector bundle~$\eta$ over a base~$B$.
We denote by $RG_k^n(\eta)$ the bundle associated to~$\eta$ with
fiber $RG_k^n$.

Let $\xi$, $\zeta$ be real vector bundles of dimensions $n_1$,~$n_2$,
respectively, over~$B$, $(n_1+n_2=n)$. Let us define the action of the positive real
group~$\Bbb R_+$ on $\xi \times \zeta$ by the formula $t\circ(v,u)=(v,tu)$.
This action is canonically extended on grassmannization
$RG^n_k(\xi\times\zeta)$. We define the vector field~$s$ by
$$
s(z)=\frac{d}{dt}\biggr|_{t=1}t\circ z
$$
where $z\in RG^n_k(\xi\times\zeta)$. Let us show that this is a
Morse--Bott vector field. It is enough to check that on the typical fiber
$RG^n_k(\Bbb R^{n_1}\times\Bbb R^{n_2})$ of the bundle
$RG^n_k(\xi\times\zeta)$ this field is a Morse-Bott field.

Zeroes of restriction of the vector field~$s$ on
$RG^n_k(\Bbb R^{n_1}\times\Bbb R^{n_2})$
are fixed points of the action~$\Bbb R_+$ on this fiber, that is, the union of
submanifolds
$RG^{n_1}_{k_1}(\Bbb R^{n_1})\times RG^{n_2}_{k_2}(\Bbb R^{n_2})$,
\ $k_1+k_2=k$ (see~\cite{7}).

In a neighborhood of a point $z_0\in RG^{n_1}_{k_1}(\Bbb R^{n_1})\times
RG^{n_2}_{k_2}(\Bbb R^{n_2})\subset
RG^n_k(\Bbb R^{n_1}\times\Bbb R^{n_2})$
one can introduce the coordinate system
$$
A(z)=\pi^{\perp}_z\circ\pi^{-1}_z
$$
where $A(z)$ is a matrix of dimension $(n-k)\times k$; \ $\pi_z$ is the
orthogonal projection of the plane~$z$ onto the plane~$z_0$; \ $\pi^{\perp}_z$
is the orthogonal projection of the plane~$z$ on the orthogonal complement
to~$z_0$ (the neighborhood of the point~$z_0$ consists of all points~$z$ for
which $\pi_{z}$ is isomorphizm). Since the plane~$z_0$ is the direct sum
of two planes $x_0\in RG^{n_1}_{k_1}(\Bbb R^{n_1})$ and
$y_0\in RG^{n_2}_{k_2}(\Bbb R^{n_2})$, the matrix $A(z)$ is canonically
divided into the blocks
$$
A(z)=\pmatrix D_1&D_2 \\ D_3&D_4 \endpmatrix
$$
where~$D_1$ is a matrix of dimension $(n_1-k_1)\times k_1$, \ $D_2$ has
dimension $(n_1-k_1)\times k_2$, \ $D_3$ has dimension $(n_2-k_2)\times k_1$,
\ $D_4$ has dimension $(n_2-k_2)\times k_2 $. The points of the submanifold
$RG^{n_1}_{k_1}(\Bbb R^{n_1})\times RG^{n_2}_{k_2}(\Bbb R^{n_2})$
in this neighborhood
are extracted by the equations $D_2=D_3=0$, and the tubular neighborhood
fiber's
points of the same submanifold over the point~$z_0$ are extracted by the
equations $D_1=D_4=0$. In this neighborhood of the point $z_0$ the
group~$\Bbb R_+$ acts by the following way
$$
t\circ\pmatrix D_1&D_2 \\ D_3&D_4 \endpmatrix=
\pmatrix D_1&D_2/t \\ tD_3&D_4 \endpmatrix.
\tag 5.1
$$

Hence, the vector field constructed is a Morse--Bott vector field.

Let us describe behavior of the vector field~$s$ on the whole space
$RG^n_k(\xi\times \zeta)$. Let $\Cal E$ be the principal $O(n)$-bundle
associated with the vector bundle $\xi\times\zeta$ over~$B$.
Zeroes of the field~$s$ are subbundles $RG^{n_1}_{k_1}(\xi)\times
RG^{n_2}_{k_2}(\zeta)$,\  $k_1+k_2=k$. The  bundle of tangents along the
fibers of the bundle $RG^n_k(\xi\times \zeta)$ is the vector bundle
$$
{\Cal E}\times_{O(n)}\tau(RG^n_k).
$$
Let $\nu(RG^{n_1}_{k_1}\times RG^{n_2}_{k_2})$ be the normal bundle
of the embedding $RG^{n_1}_{k_1}\times RG^{n_2}_{k_2}\subset
RG^n_k$. It is isomorphic to the direct sum of the vector bundles
$\xi_{n_1-k_1}\otimes \xi_{k_2}$ and $\xi_{k_1}\otimes \xi_{n_2-k_2}$, where
we denote by~$\xi_k$ the $k$-dimensional tautological vector bundle
over $RG^n_k$; \ $\xi_{n-k}$ is the~$(n-k)$-dimensional vector bundle
over $RG^{n_2}_{k_2}$ that is complementary to~$\xi$ (remark that the
bundle~$\xi_{n-k}$ is isomorphic to the $(n-k)$-dimensional tautological vector
bundle over $RG^n_{n-k}$). Thus, the normal bundle of the embedding
$$
RG^{n_1}_{k_1}(\xi)\times RG^{n_2}_{k_2}(\zeta)\subset
RG^n_k(\xi\times\zeta)
$$
is the bundle
$$
\nu={\Cal E}\times_{O(n)}\nu(RG^{n_1}_{k_1}\times RG^{n_2}_{k_2}),
$$
and
$$
\nu\cong\pi^*_1\xi(k_1)\otimes\pi^*_2\zeta(n_2-k_2)
\oplus\pi^*_1\xi(n_1-k_1)\otimes\pi^*_2\zeta(k_2)
$$
where we denote by $\xi(k)$ the~$k$-dimensional tautological vector bundle
over $RG^n_k(\xi)$ (respectively, $\xi(n-k)$ is the $(n-k)$-dimensional
tautological vector bundle over $RG^n_{n-k}(\xi)\cong
RG^n_k(\xi)$); \ $\pi_1$,~$\pi_2$ are the projections on the first and the
second factors in the product
$RG^{n_1}_{k_1}(\xi)\times RG^{n_2}_{k_2}(\zeta)$.

\proclaim{Lemma 5.1}
The vector field~$s$ associates the map {\rm(}see~\S\,{\rm3)}
$$
\psi^s_{k_1k_2}\:\nu\to\nu
$$
to the subbundle $RG^{n_1}_{k_1}(\xi)\times RG^{n_2}_{k_2}(\zeta)$
by the formula
$$
\phi^s_{k_1k_2}(u,v)=(u,-v)
$$
where $u\in\pi^*_1\xi(k_1)\otimes\pi^*_2\zeta(n_2-k_2)$,
\ $v\in\pi^*_1\xi(n_1-k_1)\otimes\pi^*_2\zeta(k_2)$.
\endproclaim

\demo{Proof}
It is enough to verify this statement on each fiber of the
bundle $RG^{n_1}_{k_1}(\xi)\times RG^{n_2}_{k_2}(\zeta)$. Owing
to~\thetag{5.1} on each fiber of this bundle restriction of~$\phi^s_{k_1k_2}$
on the fiber
is the symmetry respectively to the subbundle $\xi_{k_1}\otimes\xi_{n_2-k_2}$.
Hence, the lemma's statement follows.
\enddemo

In the previous notation we have the following statments.

\proclaim{Corollary 5.1}
In the case $k_2=0$ the index of the vector field~$s$ on
the subbundle of zeroes $RG^n_k(\xi)\times B$ is equal to~{\rm1}.
\endproclaim

The proof is obvious.

\proclaim{Corollary 5.2}
The square of the index $\Ind_s\bigl(RG^{n_1}_{k_1}(\xi)\times
RG^{n_2}_{k_2}(\zeta)\bigr)$ of the vector field~$s$ zero is equal
to~{\rm1} for any~$k_1$ and~$k_2$ whenever $k_1+k_2=k$.
\endproclaim

\demo{Proof}
Let us remark that $(\psi^s_{k_1k_2})^2=\text{id}$ for all~$k_1$ and~$k_2$
such that $k_1+k_2=k$. Thus, the statement of the corollary follows
from the definition of the multiplication in the ring $\pi^0_S(\cdot)$
as a composition of the maps.
\enddemo

\proclaim{Corollary 5.3}
The following equality holds
$$
\epsilon\bigl(\Ind_s\bigl(RG^{n_1}_{k_1}(\xi)\times
RG^{n_2}_{k_2}(\zeta)\bigr)\bigr)=(-1)^{(n_1-k_1)k_2}.
$$
\endproclaim

The proof is obvious.

In the end of this section we shall calculate the index
of the zero of the vector field~$s$ in the case when the corresponding
map~$\psi^s$ is symmetry respectively to a certain subbundle. As a
corollary we shall get an explicit expression for zero indices
of the vector field~$s$ constructed in this section.

In~\cite{6} the cobordism theory $CO^*(\cdot)$
constructed using manifolds, whose stable tangent bundle is endowed with
a fixed structure of complexification of the real vector bundle, was
considered. Any bundle of the type $\Bbb C\otimes\xi$ is oriented in this
theory, that is, it possesses the Thom class. For any~$n$-dimensional vector
bundle~$\xi$ the canonical Thom class of the bundle $\Bbb C\otimes\xi$ in
this theory is defined using the classified map
$$
u(\xi)\:T(\Bbb C\otimes\xi)\to T(\Bbb C\otimes\xi_n)
$$
where $\xi_n$ is the $n$-dimensional tautological vector bundle over $BO(n)$.
Let us denote by $\bar u(\xi)$ another Thom class of the bundle
$\Bbb C\otimes\xi$ defined in this theory by the map
$$
\bar u(\Bbb C\otimes\xi):T(\Bbb C\otimes\xi)=T(\xi\oplus\xi)
\overset{(1,-1)}\to{\longrightarrow}T(\xi\oplus\xi)=
T(\Bbb C\otimes\xi)\overset{u(\xi)}\to{\longrightarrow}T(\Bbb C\otimes\xi_n).
$$
Let $\xi$ be a real vector bundle over a finite cell complex~$B$ and
$\xi^{\perp}$ be its orthogonal complement. Let us construct the map
$$
\rho(\xi)\:T(\xi\oplus\xi^{\perp})\overset{(-1,1)}\to{\longrightarrow}
T(\xi\oplus\xi^{\perp}).
$$
We denote by $\sigma_m$ the suspension isomorphism
$CO^*(B^+)\to CO^{*+m}(S^m\wedge B^+)$.

\proclaim{Theorem 5.1}
For the real vector bundle~$\xi$ in the ring $CO^*(B^+)$
the following formula holds
$$
\sigma^{-1}_k\rho(\xi)^*\sigma_k(1)=u(\xi)^{-1}\bar u(\xi)
\tag 5.2
$$
where we denote by $u^{-1}$ the inverse map to the Thom isomorphism
in the theory $CO^*(\cdot)$ induced by the Thom class~$u$, where
$k=\dima(\xi\oplus\xi^{\perp})$.
\endproclaim

\demo{Proof}
Remark that the map $\Sigma^k\bigl(\rho(\xi)\bigr)$is homotopic
to the map
$$
T(\xi\oplus\xi\oplus\xi^{\perp}\oplus\xi^{\perp})
@>(1,-1,1,1)>>
T(\xi\oplus\xi\oplus\xi^{\perp}\oplus\xi^{\perp}).
$$
Thus, the statement of the theorem immediately follows from multiplicativity
of the Thom class.
\enddemo

We shall denote below by
$\gamma(\xi)$ the element~\thetag{5.2} of the ring $CO^*(B^+)$.

Suppouse~$F_l$ is a connected component of the zero submanifold of the
Morse--Bott field~$s$ on the manifold~$F$, and the corresponding map~$\psi^s$
is a symmetry respectively to a certain subbundle. Now let $\xi$ be the
subbundle on which~$\psi^s$ acts by fiber multiplication by~$-1$.

\proclaim{Corollary 5.4}
The following formula holds
$$
\bigl(\Ind_s(F_l)\bigr)_{CO}=\gamma(\xi).
$$
\endproclaim

In particular, for the vector~field $s$ on grassmannization of the splittable
vector bundle constructed above in this section the equality holds
$$
\bigl(\Ind\bigl(RG^{n_1}_{k_1}(\xi)\times
RG^{n_2}_{k_2}(\zeta)\bigr)\bigr)_{CO}=
\gamma\bigl(\pi^*_1\bigl(\xi(n_1-k_1)\bigr)\otimes
\pi^*_2\bigl(\zeta(k_2)\bigr)\bigr).
$$

Let us consider the formal group
$f(u,v)=u+v+\sum_{i,j\geq1}\alpha_{ij}u^iv^j$, where
$\alpha_{ij}\in \Omega_U$ in complex cobordisms~\cite{15}.
Let $\bar u\in U^2\bigl(CP(\infty)\bigr)$ be the inverse element for
$u\in U^2\bigl(CP(\infty)\bigr)$ in this group, that is, such a formal row in
$\Omega_U\bigl[[u]\bigr]$ that $f(u,\bar u)=0$. Denote by $\phi(x)$ the row
such that $\bar u=u\phi(u)$, and $\phi_n(c_1,\dots,c_n)$ the row in the ring
$U^*\bigl(BU(n)\bigr)\cong \Omega_U\bigl[[c_1,\dots,c_1]\bigr]$
defined by the identity
$$
\phi_n\bigl(\sigma_1(x_1,\dots,x_n),\dots,\sigma_n(x_1,\dots,x_n)\bigr)=
\prod^n_{i=1}\phi(x_i)
$$
where $\sigma_i$ are elementary symmetric polynomials of~$n$ variables.

Using the canonical transformation
$\mu^{CO}_U\:CO^*(\cdot)\to U^*(\cdot)$, immediately from Theorem~5.1 we get.

\proclaim{Corollary 5.5}
In the ring of complex cobordisms $U^*(B^+)$ the
following formula holds
$$
\bigl(\Ind_s(F_l)\bigr)_U=\gamma_U(\xi)=
\phi_n\bigl(c_1(\Bbb C\xi),\dots,c_n(\Bbb C\xi)\bigr)
$$
where $\gamma_U(\xi)=\mu^{CO}_U\bigl(\gamma(\xi)\bigr)$.
\endproclaim

\proclaim{Corollary 5.6}
For a one-dimensional real vector bundle~$\xi$
in the complex cobordism ring $U^*(B^+)$ the following formula
holds
$$
\gamma_U(\xi)=-1-\sum_{i,j\geq1}\alpha_{ij}\bigl(c_1(\Bbb C\xi)\bigr)^{i+j-1}.
$$
\endproclaim

\demo{Proof}
Remark that $c_1(\Bbb C\xi)=c_1(\,\overline{\Bbb C\xi}\,)$. Thus,
$$
f\bigl(c_1(\Bbb C\xi),c_1(\Bbb C\xi)\bigr)=0.
$$
Hence, the desired result is obtained.
\enddemo

\head
\S\,6. Addition Theorems for Pontrjagin Characteristic Classes
\endhead

 In this section we shall apply the localization property of the transfer
map for studying characteristic classes of the real vector bundles
introduced in~\cite{6}. In case of ordinary integral
cogomologies integer classes coincide with classical
Pontrjagin characteristic classes and that is why they were called
in~\cite{6} as Pontrjagin classes. The construction of characteristic
classes \cite{6} is based on the following properties of the
transfer map (see~\cite{3},~\cite{8}).

Let $(E_i,F_i,B_i,p)$, \ $i=1,2$, be fiber bundles whose fibers are
smooth compact closed manifolds.

\proclaim{Property 6.1 \rm(functoriality)}
If
$$
\psi\:(E_1,F_1,B_1,p)\to(E_2,F_2,B_2,p)
$$
is a bundle map, then
$\tau(p_2)\circ\psi\sim\psi\circ\tau(p_1)$.
\endproclaim

\proclaim{Property 6.2 \rm(multiplicativity)}
If
$$
(E,F,B,p)=(E_1\times E_2,F_1\times F_2,B_1\times B_2,p_1\times p_2),
$$
then $\tau(p_1)\wedge\tau(p_2)\sim\tau(p_1\times p_2)$.
\endproclaim

Let $\xi$ be a real~$n$-dimensional vector bundle over a finite cell
complex~$B$ (we can assume without loss of generality that $B$ is a smooth
manifold).  Let us consider grassmannization
of the vector bundle~$\xi$ , that is, the bundle
$$
\bigl(RG^n_k(\xi),RG^n_k,B,p_k\bigr).
$$
for given~$k$, \ $1\le k\le n$.
Let $\tau(p_k)$ be the transfer map for this bundle.
We denote by $\xi(k)$ the tautological vector bundle over $RG^n_k(\xi)$ and by
$f_k$ its classified map
$$
\xi(k)=f^*_k \xi_k
$$
where $\xi_k$ is the tautological bundle over $BO(k)$. Let
$$
s_0\:BO(k)\to T(\Bbb C\otimes\xi_k)
$$
be the embedding on the zero section.

\definition{Definition 6.1}
In the theory $CO^*$ the half-integer Pontrjagin class $p_{k/2}(\xi)$
is equal to
$$
p_{k/2}(\xi)=\tau(p_k)^*f^*\chi(\Bbb C\otimes\xi_k)=\tau(p_k)^*f^*[s_0]
\tag 6.1
$$
where $\chi(\Bbb C\otimes\zeta)$ is the Euler class of the bundle
$\Bbb C\otimes\zeta$ in the theory $CO^*$ (see~\cite{6}).
Suppose $p_{k/2}(\xi)=0$ for $k > n$.
\enddefinition

As it follows from the property of the transfer map~6.1 the
formula~\thetag{6.1} defines characteristic classes of real vector bundles.

\proclaim{Theorem 6.1}
Let $\xi$~and $\zeta$ be real vector bundles
{\rm($\dima\xi=n_1$, \ $\dima\zeta=n_2$, $n_1+n_2=n$)}.
Then  the following formula holds
$$
\split
p_{k/2}(\xi\times\zeta)&=\sum_{k_1+k_2=k}\tau(p_{k_1k_2})^*
\bigl(\gamma\bigl(\pi^*_1\xi(n_1-k_1)\otimes\pi^*_2\zeta(k_2)\bigr)
\chi\bigl(\Bbb C\otimes\xi(k_1)\bigr)
\\
&\qquad\times
\chi\bigl(\Bbb C\otimes\zeta(k_2)\bigr)\bigr)
\endsplit
$$
for any~$k$,
where $\tau(p_{k_1k_2})$ is the transfer map of the bundle
$RG^{n_1}_{k_1}(\xi)\times RG^{n_2}_{k_2}(\zeta)$, \ $\pi_1$ and $\pi_2$
are the projections on the first and second factors in the product
$RG^{n_1}_{k_1}(\xi)\times RG^{n_2}_{k_2}(\zeta)$.
\endproclaim

\demo{Proof}
Let us apply the localization property of the transfer map
in the form of Theorem~4.2 to the vector field~$s$ on
$RG^n_k(\xi\times\zeta)$ constructed in~\S\,5. Owing to
Corollary~5.4 we obtain
$$
\split
\tau(p_k)^*\chi\bigl(\Bbb C\otimes\eta(k)\bigr)&=
\sum_{k_1+k_2=k}\tau(p_{k_1k_2})^*
\bigl(\gamma\bigl(\pi^*_1\xi(n_1-k_1)\otimes\pi^*_2\zeta(k_2)\bigr)
\chi\bigl(\Bbb C\otimes\xi(k_1)\bigr)
\\
&\qquad\times
\chi\bigl(\Bbb C\otimes\zeta(k_2)\bigr)\bigr)
\endsplit
\tag 6.2
$$
where we denote by $\eta(k)$ the~$k$-dimensional tautological vector bundle
over $RG^n_k(\xi\times\zeta)$. It remains to remark that the left part
of~\thetag{6.2} is equal to $p_{k/2}(\xi\times\zeta)$ by the definition.
\enddemo

Let us show the method to obtain from Theorem 6.1 the addition formula
for Pontrjagin classes modulo elements of 2-primary order.

\proclaim{Theorem 6.2}
Let $\xi$ and $\zeta$ be real vector bundles
over a finite cell
base~$B$ \ $(\dima\xi=n_1$, \ $\dima\zeta=n_2$, \ $n_1+n_2=n)$.
For any~$k$ the following formula holds
$$
p_{k/2}(\xi\oplus\zeta)=\sum_{k_1+k_2=k}
(-1)^{(n_1-k_1)k_2}p_{k_1/2}(\xi)p_{k_2/2}(\zeta)
\pmod{2\Tors}.
$$
\endproclaim

\demo{Proof}
Owing to functoriality of the Pontrjagin classes
it is enough to prove the statement of the theorem for the bundle
$\pi^*_1(\xi)\oplus\pi^*_2(\zeta)$ over $B\times B$, where $\pi_1$~and
$\pi_2$ are the projections on the first and the second factors in the
product $B\times B$. Remark also that
$$
\pi^*_1(\xi)\oplus\pi^*_2(\zeta)\cong\xi\times\zeta.
$$
Assume that
$$
\alpha=\bigl(\Ind_s\bigl( RG^{n_1}_{k_1}(\xi)\times
 RG^{n_2}_{k_2}(\zeta)\bigr)\bigr)_{CO}.
$$
>From Corollary~5.2 \ $\alpha^2=1$. We have $\alpha=\pm 1+u$, where
$u\in\widetilde{CO^0}\bigl(RG^{n_1}_{k_1}(\xi)\times
RG^{n_2}_{k_2}(\zeta)^+\bigr)$.
Since we work in the category of finite cell complexes, the element~$u$
is nilpotent. The relation $\alpha^2=1$ involves the equality $\pm 2u+u^2=0$.
Hence and from nilpotency of~$u$ we conclude that the element~$u$ is
2-primary, that is, $u\in 2\Tors$. Thus, $\alpha=\pm 1(\bmod(2\Tors))$. From
Corollary~5.3
$$
\alpha\equiv(-1)^{(n_1-k_1)k_2}\pmod{2\Tors}.
\tag 6.3
$$

We denote by $\eta(k)$ the $k$-dimensional tautological vector bundle over
$RG^n_k(\xi\times\zeta)$. From the formula~\thetag{6.2} taking into
account~\thetag{6.3} we obtain
$$
\split
&\tau(p_k)^*\chi\bigl(\Bbb C\otimes\eta(k)\bigr)=
\\
&\qquad=
\sum_{k_1+k_2=k}(-1)^{(n_1-k_1)k_2}\tau(p_{k_1k_2})^*
\bigl(\chi(\Bbb C\otimes\xi_{k_1})\times\chi(\Bbb C\otimes\zeta_{k_2})\bigr)
\pmod{2\Tors}.
\endsplit
\tag 6.4
$$

>From multiplicativity of the transfer map 6.2 it follows that
$$
\split
&\tau(p_{k_1}\times p_{k_2})^*
\bigl(\chi\bigl(\Bbb C\otimes\xi(k_1)\bigr)\times
\chi\bigl(\Bbb C\otimes\zeta(k_2)\bigr)\bigr)
\\
&\qquad=
\tau(p_{k_1})^*\chi\bigl(\Bbb C\otimes\xi(k_1)\bigr)
\tau(p_{k_2})^*\chi\bigl(\Bbb C\otimes\zeta(k_2)\bigr).
\endsplit
\tag 6.5
$$

Combining~\thetag{6.4} and~\thetag{6.5} we obtain the formula required
$$
p_{k/2}(\xi\times\zeta)=\sum_{k_1+k_2=k}(-1)^{(n_1-k_1)k_2}
p_{k_1/2}(\xi)\times p_{k_2/2}(\zeta)\pmod{2\Tors}.
$$
\enddemo

\proclaim{Theorem 6.3}
Half-integer Pontrjagin classes are stable
$$
p_{k/2}(\xi\oplus1)=p_{k/2}(\xi)=p_{k/2}(1\oplus\xi).
$$
\endproclaim

\demo{Proof}
{\Bad
Let us apply the localization property of the transfer map
to the bundle $RG^{n+1}_k(\xi\oplus 1)$ \ $(\dim\xi=n)$. Using
Corollary~5.1 we obtain that for any element
$x\in CO^*\bigl(RG^{n+1}_k(\xi\oplus1)^+\bigr)$ and for a certain element
$\alpha\in
CO^0\bigl(RG^{n+1}_k(\xi\oplus1)^+\bigr)$ the following splitting holds
}
$$
\tau(p)^*x=\tau(p_1)^*i^*_1x+\tau(p_2)^*\alpha i^*_2x
$$
where $\tau(p_1)$ is the transfer map of the bundle
$\bigl(RG^n_k(\xi),RG^n_k,B,p_1\bigr)$ and $\tau(p_2)$
is the transfer map of the bundle
$\bigl(RG^n_{k-1}(\xi),RG^n_{k-1},B,p_2\bigr)$,
the other maps were defined in the formulation of
the localization property.

Since under the embedding
$RG^n_{k-1}(\xi)\subset RG^{n+1}_k(\xi\oplus1)$ the tautological vector
bundle $\xi(k)$ over $RG^{n+1}_k(\xi\oplus 1)$ turns to the sum
$\xi(k-1)\oplus 1$, where $\xi(k-1)$ is the tautological vector bundle
over $RG^n_{k-1}(\xi)$, then
$$
\split
\tau(p)^*\chi\bigl(ó\otimes\xi(k)\bigr)&=
\tau(p_1)^*\chi\bigl(ó\otimes\xi(k)\bigr)+
\tau(p_2)^*\alpha\chi\bigl(ó\otimes\bigl(\xi(k-1)\oplus1\bigr)\bigr)
\\
&=
\tau(p_1)^*\chi\bigl(ó\otimes\xi(k)\bigr).
\endsplit
$$
Hence, we obtain the first equality of the theorem.

To prove the second
equality it is enough to multiply the vector field from~\S\,5 by~$-1$
and to proceed with an analogous reasoning.
\enddemo

\proclaim{Corollary 6.1}
 The half-integer Pontrjagin class $p_{k/2}(\xi)$ \
is {\rm2}-primary for odd~$k$.
\endproclaim

\demo{Proof}
Using Theorem~6.2
$$
p_{k/2}(1\oplus\xi)=(-1)^kp_{k/2}(\xi)\pmod{2\Tors},
$$
and from Theorem~6.3
$$
p_{k/2}(1\oplus\xi)=p_{k/2}(\xi).
$$
For odd~$k$ these equalities are compatible only if
$p_{k/2}(\xi)\in2\Tors$.
\enddemo

\proclaim{Corollary 6.2}
For any two real vector bundles~$\xi$ and~$\zeta$
the following relation holds~\cite{6}
$$
p_{2k/2}(\xi\oplus\zeta)=\sum_{k_1+k_2=k}p_{2k_1/2}(\xi)
p_{2k_2/2}(\zeta)\pmod{2\Tors}.
$$
\endproclaim

\proclaim{Theorem 6.4}
For the image of Pontrjagin characteristic classes in
the complex cobordisms there is a formula, which expresses the
half-integer Pontrjagin classes of the sum of two real vector bundles
in terms of characteristic classes of their summands.
\endproclaim

\demo{Proof}
Owing to splitting~\thetag{6.2} it is enough to prove
that $\bigl(\Ind_s\bigl( RG^{n_1}_{k_1}(\xi)\times
RG^{n_2}_{k_2}(\zeta)\bigr)\bigr)_U$
can be expressed as a polynomial of characteristic classes of the bundles
$\xi(k_1)$ and $\zeta(k_2)$ with coefficients in
$p^*_1\bigl(U^*(B)\bigr)\otimes p^*_2\bigl(U^*(B)\bigr)$, where~$p_1$ and
$p_2$ are the projections
in the bundles $RG^{n_1}_{k_1}(\xi)$ and $RG^{n_2}_{k_2}(\zeta)$ respectively.
Owing to Lemma~5.1 and Corollary~5.5
the following equality holds in complex cobordisms
$$
\bigl(\Ind_s\bigl( RG^{n_1}_{k_1}(\xi)\times
 RG^{n_2}_{k_2}(\zeta)\bigr)\bigr)_U=
\phi_{(n_1-k_1)k_2}\bigl(c_1(\eta),\dots,c_{(n_1-k_1)k_2}(\eta)\bigr)
$$
where $\eta\cong \Bbb C\otimes\bigl(\xi(n_1-k_1)\otimes\zeta(k_2)\bigr)$.
Using the splitting principle for complex vector bundles and the formula for
the formal group
$$
f\bigl(c_1(\eta_1),c_1(\eta_2)\bigr)=c_1(\eta_1\otimes\eta_2)
$$
we obtain the expression of the classes $c_i(\eta)$ in terms of characteristic
classes of the bundles $\Bbb C\otimes\xi(n_1-k_1)$ and
$\Bbb C\otimes\zeta(k_2)$. Finally, it is enough to remark that
$\xi(n_1-k_1)\oplus\xi(k_1)\cong p^*_1\xi$.
\enddemo

\remark{Remark 6.1}
Half-integer Pontrjagin classes $p_{k/2}(\xi)$ with
odd numbers~$k$ can have any 2-primary order. For example, let $\xi_1$ be the
one-dimensional tautological vector bundle over $RP^{2n+1}$. We denote
by $\mu^{CO}_K\:CO^*(\cdot)\to K^*(\cdot)$ the canonical transformation of
the cogomology theories. The order of
$\mu^{CO}_K\bigl(p_{1/2}(\xi_1)\bigr)\in K^*( RP^{2n+1})$ is equal to~$2^n$.
\endremark

At the end of this section we shall consider in detail the addition theorem
for the first half-integer Pontrjagin class (see~\cite{9}).

We assume that $\xi$ and $\zeta$ are vector bundles over a finite cell
base~$B$. Let
$\tau(p)$ be the transfer map of the bundle
$RP(\xi\oplus\zeta)$,\ $\tau(p_1)$ the transfer map of the bundle
$RP(\xi)$, \ $\tau(p_2)$ the
transfer map of the bundle $RP(\zeta)$. Let us denote by~$i_1$ the embedding
$RP(\xi)\subset RP(\xi\oplus\zeta)$, by $i_2$ the embedding
$RP(\zeta)\subset RP(\xi\oplus\zeta)$.

\proclaim{Theorem 6.5}
 The following
formula holds
$$
\tau(p)^*(x)=\tau(p_1)^*i^*_1(x)+\tau(p_2)^*
\bigl(\gamma_U\bigl(p^*_2\xi\otimes\zeta(1)\bigr)i^*_2(x)\bigr)
$$
for any $x\in U^*\bigl(RP(\xi\oplus\zeta)\bigr)$.
\endproclaim

\demo{Proof}
The statement of the theorem immediately follows from
Theorem~4.2 and Corollaries~5.1 and~5.4.
\enddemo

\proclaim{Corollary 6.3}
For any two one-dimensional real vector bundles~$\xi$
and~$\zeta$ over a finite cell base the following equality
holds
$$
p_{1/2}(\xi\oplus\zeta)=u-v-\sum_{i,j\geq1}\alpha_{ij}f(u,v)^{i+j-1}v
$$
where $u=p_{1/2}(\xi)$, \ $v=p_{1/2}(\zeta)$.
\endproclaim

\demo{Proof}
For any one-dimensional vector bundle $\xi$ we have
$p_{1/2}(\xi)=c_1(\Bbb C\otimes\xi)$. Thus, the statement required follows
from the explicit form $\gamma_U(\xi\otimes\zeta)$ (see Corollary~5.6).
\enddemo

Let $[u]_2=ua(u)$, where
$a(u)=2+\sum_{i,j\geq1} \alpha_{ij} u^{i+j-1}$. Then the result of
Corollary~6.3 can be written in the form
$$
p_{1/2}(\xi\oplus\zeta)=u+v-a\bigl(f(u,v)\bigr)v.
$$
Next let us introduce the following series, symmetric on~$u$ and~$v$:
$$
\delta(u,v)=\frac{a(u)-a(v)}{u-v},\qquad
d(u,v)=\frac{va(u)-ua(v)}{u-v}.
$$
Let us write the series $\alpha(u)=1+\sum_{i\geq1}\alpha_{i1}u^i$
in the form
$$
\alpha(u)=\alpha_0(u^2)+u\alpha_1(u^2).
$$

\proclaim{Theorem 6.6}
For any two real one-dimensional vector bundles~$\xi$
and~$\zeta$ the following equality holds
$$
p_{1/2}(\xi\oplus\zeta)=
u+v-uv\bigl[\alpha_0(uv)\delta(u,v)+\alpha_1(uv)d(u,v)\bigr]
$$
where $u=p_{1/2}(\xi)$ and $v=p_{1/2}(\zeta)$.
\endproclaim

\demo{Proof}
Let
$$
\Phi(u,v)=1+\sum_{i,j\geq1}\alpha_{ij}u^i
\left(\frac{v^j-{\bar u}^j}{v-\bar u}\right)\in\Omega_U\bigl[[u,v]\bigr]
$$
where $\bar u$ satisfies the relation $f(u,\bar u)=0$. Then
$$
f(u,v)=f(u,v)-f(u,\bar u)=(v-\bar u)\Phi(u,v).
$$
We get
$$
f\bigl([u]_2,[v]_2\bigr)=\bigl[f(u,v)\bigr]_2=f(u,v)a\bigl(f(u,v)\bigr).
$$
Therefore,
$$
\bigl([v]_2-[\bar u]_2\bigr)\Phi\bigl([u]_2,[v]_2\bigr)=
(v-\bar u)\Phi(u,v)a\bigl(f(u,v)\bigr).
$$
Thus, in the ring $\Omega_U\bigl[[u,v]\bigr]$
the following identity holds
$$
\Phi(u,v)a\bigl(f(u,v)\bigr)v=
\left([v]_2+\bar uv\frac{a(v)-a(\bar u)}{v-\bar u}\right)
\Phi\bigl([u]_2,[v]_2\bigr).
$$

To prove Theorem 6.6 it is enough to assume that $\xi$ and $\zeta$ are
one-dimensional tautological vector bundles over $RP(\infty)$. Thus,
all further computations can be held in the factor ring
$$
A=\Omega_U\bigl[[u,v]\bigr]/\bigl([u]_2,[v]_2\bigr).
$$
>From the identity $[u]_2=(u-\bar u)\Phi(u,u)$ in $\Omega_U\bigl[[u]\bigr]$
immediately follows that in the factor ring $u=\bar u$ and $v=\bar v$.
Using the relations $[u]_2=0$ and $u=\bar u$ we obtain the following
identity in~$A$
$$
\Phi(u,v)a\bigl(f(u,v)\bigr)v=uv\,\delta(u,v),
$$
from which it immediately follows that in~$A$ the following relation holds
$$
a\bigl(f(u,v)\bigr)v=a\bigl(f(u,v)\bigr)u.
$$

Using now an explicit form of the series $\Phi(u,v)$ introduced above we
obtain that in the ring~$A$ the following identity holds
$$
a\bigl(f(u,v)\bigr)\Phi(u,v)=
a\bigl(f(u,v)\bigr)\frac{\partial}{\partial v}f(u,v).
$$
\enddemo

\proclaim{Lemma 6.1}
In the ring $\Omega_U\bigl[[u,v]\bigr]$ the following formula holds
$$
\frac{\partial}{\partial v}f(u,v)=\frac{CP(v)}{CP\bigl(f(u,v)\bigr)}
$$
where $CP(v)={dg(v)}/{dv}=1+\sum_{n\geq1}\bigl[CP(n)\bigr]v^n$ and
$g(v)$ is logarithm of the group $f(u,v)$.
\endproclaim

\demo{Proof}
In the ring $\Omega_U\otimes Q\bigl[[u,v]\bigr]$ for the series
$f(u,v)$ the following expression holds
$$
f(u,v)=g^{-1}\bigl(g(u)+g(v)\bigr).
$$
Since
$$
\frac{dg^{-1}(x)}{dx}\cdot\frac{dg\bigl(g^{-1}(x)\bigr)}{dg^{-1}(x)}=1,
$$
we obtain that
$$
\frac{\partial}{\partial v}f(u,v)=
\frac{\partial g^{-1}\bigl(g(u)+g(v)\bigr)}
{\partial g(v)}\frac{\partial g(v)}{\partial v}=
\frac{CP(v)}{CP\bigl(f(u,v)\bigr)}.
$$

Returning to the ring~$A$ we obtain
$$
a\bigl(f(u,v)\bigr)\frac{\partial}{\partial v}f(u,v)=
a\bigl(f(u,v)\bigr)\frac{CP(v)}{CP\bigl(f(u,v)\bigr)}=
a\bigl(f(u,v)\bigr)CP(v).
$$
Let us remind that
$$
\frac{1}{CP(v)}=
\frac{\partial f(u,v)}{\partial u}\biggr|_{u=0}=
1+\sum_{i\geq1}\alpha_{i1}v^i=\alpha(v).
$$
Therefore,
$$
uv\,\delta(u,v)\alpha(v)=
\alpha(v)\Phi(u,v)a\bigl(f(u,v)\bigr)v=a\bigl(f(u,v)\bigr)v.
$$
\enddemo

\proclaim{Lemma 6.2}
In the ring~$A$ the following relations hold
$$
\align
uv\,\delta(u,v)v&=uv\,d(u,v),
\\
uv\,\delta(u,v)v^2&=uv\,\delta(u,v)uv.
\endalign
$$
\endproclaim

\demo{Proof}
The first equality follows from the relation
$$
\split
uv\frac{a(v)-a(u)}{v-u}v&=uv\frac{va(v)-va(u)+ua(v)-ua(v)}{v-u}
\\
&=uv\frac{ua(v)-va(u)+(v-u)a(v)}{v-u}=
uv\frac{ua(v)-va(u)}{v-u}=uv\,d(u,v).
\endsplit
$$
The second one follows from the relation
$a\bigl(f(u,v)\bigr)v=a\bigl(f(u,v)\bigr)u$:
$$
uv\,\delta(u,v)v^2=\Phi(u,v)a\bigl(f(u,v)\bigr)v^3=
\Phi(u,v)a\bigl(f(u,v)\bigr)vuv=uv\,\delta(u,v).
$$

Applying the relations from Lemma 6.2 we obtain
$$
\multline
uv\,\delta(u,v)\alpha(v)=
uv\,\delta(u,v)\bigl(\alpha_0(v^2)+v\alpha_1(v^2)\bigr)
\\
=
uv\bigl(\alpha_0(uv)\delta(u,v)+\alpha_1(uv)\delta(u,v)v\bigr)=
uv\bigl(\alpha_0(uv)\delta(u,v)+\alpha_1(uv)d(u,v)\bigr),
\endmultline
$$
that allows us to finish the proof of Theorem~6.6.
\enddemo

This theorem gives an explicit form of the formal
series, symmetric on~$u$ and~$v$
$$
b(u,v)=u+v+\sum_{i,j\geq1}\beta_{ij}u^iv^j
$$
such that for any two one-dimensional real vector bundles~$\xi$ and~$\zeta$
over~$X$ the following formula holds
$$
p_{1/2}(\xi\oplus\zeta)=b(u,v)\in U^2(X)
$$
where $u=p_{1/2}(\xi)$ and $v=p_{1/2}(\zeta)$.

\proclaim{Theorem 6.7 \rm(the addition formula for the first half-integer
Pontrjagin class)}For any two real vector bundles~$\xi$ and~$\zeta$
$$
p_{1/2}(\xi\oplus\zeta)=
u+v+\sum_{k,l\geq1}\beta_{kl}s_{k-1}(u)s_{l-1}(v)
$$
where $u=p_{1/2}(\xi)$,\ $v=p_{1/2}(\zeta)$, \ $s_k$ are the
Landweber--Novikov operations in complex cobordisms~\cite{15}.
\endproclaim

\demo{Proof}
>From the localization theorem for the transfer map it follows
that
$$
p_{1/2}(\xi\oplus\zeta)=p_{1/2}(\xi)+
\tau(p)^*\bigl(\gamma_U\bigl(p^*(\xi)\otimes\zeta(1)\bigr)\cdot w\bigr)
$$
where $\tau(p)$ is the transfer map of the bundle $RP(\zeta)$,
\ $w=\chi\bigl(\Bbb C\otimes\zeta(1)\bigr)$.  First let us assume that
$\dima\xi=1$. We have
$$
\align
p_{1/2}(\xi\oplus\zeta)&=
u+\tau(p)^*\bigl(\gamma_U\bigl(p^*\xi\otimes\zeta(1)\bigr)w\bigr),
\tag 6.6
\\
p_{1/2}(\zeta\oplus\xi)&=v+\gamma_U(\zeta\otimes\xi)u
\tag 6.7
\endalign
$$
where $u=p_{1/2}(\xi)$, \ $v=p_{1/2}(\zeta)$.

Remark that from Corollary~5.6
$$
\split
\tau(p)^*\bigl(\gamma_U\bigl(p^*\xi\otimes\zeta(1)\bigr)w\bigr)&=
-\tau(p)^*\biggl(w+\sum_{i,j\geq1}\alpha_{ij}f(p^*u,w)^{i+j-1}w\biggr)
\\
&=
\tau(p)^*\bigl(p_{1/2}\bigl(p^*\xi\oplus\zeta(1)\bigr)-
p_{1/2}(p^*\xi)\bigr).
\endsplit
\tag 6.8
$$
The Landweber-Novikov operations are stable cohomology operations, therefore,
they commute with the map $\tau(p)^*$. Owing to Theorem~6.6 we obtain
that the expression~\thetag{6.8} is equal to
$$
\tau(p)^*\biggl(w+\sum_{k,l\geq1}\beta_{kl}(p^*u)^kw^l\biggr)=
v+\sum_{k,l\geq1}\beta_{kl}u^ks_{l-1}(v).
$$
Setting~\thetag{6.6} to be equal to~\thetag{6.7} we obtain
$$
\gamma_U(\zeta\otimes\xi^*)\cdot u=
u+\sum_{k,l\geq1}\beta_{kl}u^ls_{k-1}(v).
$$

Consider now the case of an arbitrary~$\xi$. We have
$$
\multline
p_{1/2}(\xi\oplus\zeta)=
u+\tau(p)^*\bigl(\gamma_U\bigl(p^*(\xi)\otimes\zeta(1)^*\bigr)w\bigr)
\\
=
u+\tau(p)^*\biggl(w+\sum_{k,l\geq1}\beta_{kl}w^lp^*s_{k-1}(u)\biggr)=
u+v+\sum_{k,l\geq1}\beta_{kl}s_{l-1}(v)s_{k-1}(u).
\endmultline
$$
This completes the proof.
\enddemo

\proclaim{Theorem 6.8}
The characteristic class $p_{1/2}$ defines the map
$$
p_{1/2}\:KO(X)\to U^2(X),
$$
on whose image there is an associative addition operation
$$
x\oplus y=x+y+\sum_{k,l\geq1}\beta_{kl}s_{l-1}(y)s_{k-1}(x).
$$
\endproclaim

Let us pay attention that this operation is not reduced to any formal
group.

\head
\S\,7. Pontrjagin Characteristic Classes of Complex Vector Bundles
\endhead

\proclaim{Theorem 7.1}
Let $\eta$ be a $n$-dimensional complex vector bundle.
The transfer map $\tau(p)$ of the bundle $RG^{2n}_{2k+1}(\eta)$ is homotopic
to the constant map. The transfer maps $\tau(p_1)$ and $\tau(p_2)$ of the
bundles $RG^{2n}_{2k}(\eta)$ and $CG^n_k(\eta)$ are connected by the relation
$$
\tau(p_1)\sim i\circ\tau(p_2)
$$
where $i\:CG^n_k(\eta)\to RG^{2n}_{2k}(\eta)$ is the canonical embedding.
\endproclaim

\demo{Proof}
Let us consider the bundle
$RG^{2n}_k(\eta)\overset{RG^{2n}_k}\to{\longrightarrow}B$.
The group $U(1)=S^1$
acts on the space of this bundle. This action is induced by the
complex structure in the vector bundle~$\eta$.
All hyperplanes~$\alpha$ such that for any vector $v\in \alpha$ the
vector $iv$ belongs to~$\alpha$ are fixed points of the action. Since the
vectors~$v$ and $iv$ are
linearly independent over the real number field, the fixed points can appear
whenever~$k$ is even. Applying the localization property to
the vector field corresponding to the action we obtain that the transfer map
$\tau(p)$ of the bundle $RG^{2n}_{2k+1}(\eta)$ is homotopic to the constant
map (see Remark~3.4).  In case of even~$k$ the fixed points of the action
are such hyperplanes that are realifications of the $(k/2)$-dimensional
complex subspaces of the bundle~$\eta$. Applying the localization property
we obtain $\tau(p_1)\sim i\circ\tau(p_2)$.
\enddemo

\proclaim{Corollary 7.1}
If the vector bundle~$\eta$ is stably isomorphic to
the complex vector bundle, then $ p_{(2k+1)/2}(\eta)=0$.
\endproclaim

\demo{Proof}
Let the bundle  $\eta\oplus [N]$ admit the structure of the
complex vector bundle. Then
$$
p_{(2k+1)/2}(\eta)=p_{(2k+1)/2}(\eta\oplus [N])=0.
$$

For the theory of complex self-conjugate cobordisms $SC^*(\cdot)$ the
following canonical transformations of the cogomology theories exist
$$
\mu^{CO}_{SC}\:CO^*(\cdot)\to SC^*(\cdot)\qquad\text{and}\qquad
\mu^{Sp}_{SC}\:Sp^*(\cdot)\to SC^*(\cdot)
$$
where $Sp^*(\cdot)$ is the symplectic cobordism theory
(see~\cite{6}).
\enddemo

\proclaim{Corollary 7.2}
If the vector bundle~$\eta$ is stably isomorphic
to the complex vector bundle, the image $p_{2k/2}(\eta)$ in $SC^*(B)$
under the transformation $\mu^{CO}_{Sp}$ belongs to
$\Ima\mu^{Sp}_{SC}\bigl(Sp^*(B)\bigr)\subset SC^*(B)$.
\endproclaim

\demo{Proof}
Since Pontrjagin classes are stable, without loss of
generality we can assume that the bundle~$\eta$ admits a complex structure.
Let $\tau(p_1)$ be the transfer map of the bundle $RG^{2n}_{2k}(\eta)$ and
$\tau(p_2)$ the transfer map of the bundle $CG^n_k(\eta)$.
Denote by $i\: CG^n_k\to RG^{2n}_{2k}$ the canonical embedding. From
Theorem~7.1 we obtain
$$
\split
\mu^{CO}_{SC}\bigl(p_{2k/2}(\eta)\bigr)&=
\mu^{CO}_{SC}\tau(p_1)^*\chi\bigl(\Bbb C\otimes\xi(2k)\bigr)
\\
&=
\mu^{CO}_{SC}\tau(p_2)^*i^*\chi\bigl(\Bbb C\otimes\xi(2k)\bigr)=
\mu^{CO}_{SC}\tau(p_2)^*\chi\bigl(\Bbb C\otimes\eta(k)\bigr)
\\
&=
\tau(p_2)^*\mu^{CO}_{SC}\chi\bigl(\Bbb C\otimes\eta(k)\bigr)=
\tau(p_2)^*\chi\bigl(\eta(k)\oplus\overline{\eta(k)}\,\bigr)
\endsplit
\tag 7.1
$$
where $\xi(2k)$ is the $2k$-dimensional tautological vector bundle over
$RG^{2n}_{2k}(\eta)$, and $\eta(k)$ is the $k$-dimensional tautological
vector bundle over $CG^n_k(\eta)$. Remark that for any complex vector
bundle~$\eta$ the bundle $\eta\oplus\bar\eta$ admits the canonical symplectic
structure. Therefore, the right side of~\thetag{7.1} belongs to
$\Ima\mu^{Sp}_{SC}\bigl(Sp^*(B)\bigr)$.
\enddemo


\enddocument